\documentclass[11pt]{amsart}
\usepackage{amsmath}
\usepackage{amsfonts}
\usepackage{amsthm}
\numberwithin{equation}{section}

\newtheorem{thm}{Theorem}
\newtheorem{prop}[thm]{Proposition}
\newtheorem{lem}[thm]{Lemma}

\begin{document}

\title[Convergence Radii for Eigenvalues]{Convergence
Radii for Eigenvalues of Tri--diagonal Matrices}

\author{J. Adduci}

\address{Department of Mathematics,
The Ohio State University,
 231 West 18th Ave,
Columbus, OH 43210, USA}

\email{adducij@math.ohio-state.edu}

\author{P. Djakov}
\address{Sabanci University, Orhanli,
34956 Tuzla, Istanbul, Turkey}

 \email{djakov@sabanciuniv.edu}

\author{B. Mityagin}

\address{Department of Mathematics,
The Ohio State University,
 231 West 18th Ave,
Columbus, OH 43210, USA}

\email{mityagin.1@osu.edu}

\thanks{B. Mityagin acknowledges the support of the Scientific and Technological
Research Council of Turkey and  the hospitality of Sabanci
University, April--June, 2008}


\begin{abstract}
Consider a family of infinite tri--diagonal matrices of the form $L+
zB,$ where the matrix $L$ is diagonal with entries $L_{kk}= k^2,$ and
the matrix $B$ is off--diagonal, with nonzero entries
$B_{k,{k+1}}=B_{{k+1},k}= k^\alpha, \; 0 \leq \alpha < 2.$ The
spectrum of $L+ zB$ is discrete. For small $|z|$ the $n$-th
eigenvalue $E_n (z),\; E_n (0) = n^2,$ is a well--defined analytic
function. Let $R_n$ be the convergence radius of its Taylor's series
about $z= 0.$  It is proved that
$$ R_n \leq C(\alpha) n^{2-\alpha} \quad \text{if} \; \; 0 \leq \alpha
<11/6.$$
\end{abstract}


\maketitle

\section{Introduction}

Since the famous 1969 paper of C. Bender and T. Wu \cite{BW},
branching points and  the crossings of energy levels have been
studied intensively in the mathematical and physical literature
(e.g., \cite{HS,A,DP,B01} and the bibliography there). In this paper
our goal is to analyze -- mostly along the lines of J. Meixner and F.
Sch\"afke approach \cite{MS54} -- a toy model of tri--diagonal
matrices.

We consider the operator family $L+zB,$ where $L$ and $B$ are
infinite matrices of the form
\begin{align}
\label{1.1} L=
\begin{bmatrix}
  q_1 & 0    & 0   & 0   &\cdot \\
  0   & q_2  & 0   & 0    &\cdot\\
  0   & 0    & q_3 & 0   &\cdot\\
  0   & 0    & 0   &q_4    &\cdot\\
  \cdot   & \cdot    & \cdot   &\cdot   &\cdot  \\
\end{bmatrix},
\qquad  B=
\begin{bmatrix}
  0          & b_1        & 0       & 0       &\cdot       \\
  c_1        &0           &b_2      & 0       &\cdot        \\
  0          &c_2         &0        & b_3     &\cdot        \\
  0          &0           &c_3      &0        &\cdot       \\
  \cdot          &\cdot           &\cdot        &\cdot      &\cdot        \\
\end{bmatrix}
\end{align}
with
\begin{align}
\label{1.2}
q_k &= k^2, \\
\label{1.3}
|b_k|, |c_k|  &\leq M k^{\alpha}, \\
\label{1.4} \alpha &< 2.
\end{align}
Sometimes we impose a symmetry condition:
\begin{equation}
\label{1.5}  b_k = \bar{c_k}.
\end{equation}

Under the conditions (\ref{1.2})--(\ref{1.4}) the spectrum of $L+zB$
is discrete. If $\alpha < 1$ then a standard use of perturbation
theory shows that there is $r>0$ such that for $|z|<r$
\begin{equation}
Sp(L+zB)=\{E_n(z)\}_{n=1}^{\infty}, \quad E_n (0) = n^2,
\end{equation}
where each $E_n (z)$  is well--defined analytic function in the disc
$\{z:|z| <r\}.$

If $\alpha \in [1,2),$ then in general there is no such $r>0.$ But
the fact that $n^2$ is a simple eigenvalue of $L$ guarantees (see
\cite{Kato}, Chapter 7, Sections 1-3) that for each $n$ there exists
$r_n>0$ such that, on the disc $\{z: |z|< r_n\},$ there are an
analytic function $E_n (z)$ and an analytic eigenvector function
$\varphi_n(z)$ with
\begin{align}
\label{l1}
(L+zB)\varphi_n(z) &= E_n(z) \varphi_n(z), \,\,\, |z| < r_n, \\
\varphi_n(0) &= e_n, \,\,\,E_n(0)=n^2.
\end{align}
Let
\begin{align} \label{l2} E_n(z) = \sum_{k=0}^{\infty} a_k(n)z^k
\end{align}
be the Taylor series of $E_n(z)$ about $0,$ and let $R_n,$ $0 < R_n
\leq \infty,$ be its radius of convergence. The asymptotic behavior
of the sequence $(R_n)$ is one of the main topics of the present
paper.

It may happen that $R_n >r_n. $ Then, by (\ref{l2}), $E_n(z)$ is
defined in the disc $\{z: |z| < R_n \}$ as an extension of the
analytic function (\ref{l1}) in $\{z: |z| < r_n \}.$  But are its
values $E_n(z)$  eigenvalues of $L+zB$ if $z$ is in the annulus $r_n
\leq |z| <R_n? $  The answer is positive as one can see from the
next considerations.

In a more general context let us define {\em Spectral Riemann
Surface}
\begin{align}
\label{l3} G = \{(z,E):  \;\exists g \in Dom (L), \; g \neq
0\;|\quad (L+zB)g = Eg \}.
\end{align}
This notion is justified by the following statement (coming from K.
Weierstrass, H. Poincare, T. Carlemann -- see discussions on the
related history in \cite{GoKr69,MSW80,GGK00}).

\begin{prop}
\label{prop01} If (\ref{1.1})--(\ref{1.4}) hold, then there exists a
nonzero entire function $\Phi(z,w)$ such that
\begin{align}
\label{l4} G = \{ (z,w) \in \mathbb{C}^2: \Phi(z,w)=0 \}.
\end{align}
\end{prop}

\begin{proof}
The identity
\begin{align}
\label{l5} (L+zB)g = wg, \quad g\neq 0, \quad  g \in Dom (L)
\end{align}
is equivalent to
\begin{equation}
\label{l6} (1-A(z,w))h =0 \quad  \text{with} \quad h=L^{1/2}g \in
Dom ( L^{1/2}), \;  h \neq 0,
\end{equation}
where
\begin{equation}
\label{l7} A(z,w) = -z L^{-1/2}BL^{-1/2} +wL^{-1}.
\end{equation}
Therefore, $w$ is an eigenvalue of the operator $L+zB$ if and only if
$1$ is an eigenvalue of the operator $A(z,w).$

On the space $S_1$ of trace class operators $T$ the determinant
\begin{align}
\label{l8} d(T) = \det (1-T)
\end{align}
is  well defined (see \cite{GoKr69}, Chapter 4, Section 1 or
\cite{Si79}, Chapter 3, Theorem 3.4), and $1 \in Sp (T)$ if and only
if $d (T)=0$ (see \cite{Si79}, Theorem 3.5 (b)).

Of course, the second term $L^{-1}$ in (\ref{l7}) is an operator of
trace class (even in $S_p, p>1/2$) by (1.2). But
(\ref{1.3})--(\ref{1.4}) imply that $L^{-1/2}BL^{-1/2}$ is in the
Schatten class $S_p, p>1/(2-\alpha)$; only $\alpha < 1$ would
guarantee  that it is of trace class.

However, (\ref{l8}) could be adjusted (see \cite{GoKr69} Chapter 4,
Section 2 or \cite{Si79}, Chapter 9, Lemma 9.1 and Theorem 9.2).
Namely, for any positive integer $p \geq 2$ we set
\begin{equation}
\label{l9} d_p(T) = \det(1-Q_p(T))
\end{equation}
where $$ Q_p (T) = 1-(1-T)\exp \left (T+ \frac{T^2}{2}+ \cdots +
\frac{T^{p-1}}{p-1} \right ). $$ Then $Q_p(T) \in S_1 $ if $ T \in
S_p,$ so $d_p$ is a well-defined function of $T \in S_p $ and $1 \in
Sp(T)$ if and only if $ d_p(T)=0.$

In our context we define, with $A(z,w) \in (\ref{l7})$ and
$p>1/(2-\alpha),$
\begin{align}
\label{l10} \Phi(z,w) = \det \left [( 1 - Q_p (A(z,w))  \right ].
\end{align}
Now, from Claim 8, Section 1.3, Chapter 4 in \cite{GoKr69} it
follows that $\Phi(z,w)$ is an entire function on $\mathbb{C}^2$.

The function  $\Phi $ vanishes at $(z,w)$ if and only if $1$ is an
eigenvalue of the operator $A(z,w),$ i.e., if and only if $(z,w) \in
G$. This completes the proof.
\end{proof}

In particular, the above Proposition implies that $\Phi(z,E_n(z))=0$
if $ |z| < r_n,$ so by analyticity and uniqueness $\Phi(z,E_n(z))=0 $
if $ r_n \leq |z| < R_n.$  Equivalence of the two definitions
(\ref{l3}) and (\ref{l4}) for the Spectral Riemann Surface  $G$
explains now that $E_n(z)$ is an eigenvalue function in the disc
$\{z: |z|<R_n\}$.

Our main focus in the search for an understanding of the behavior of
$R_n$ will be on the special case where
\begin{align}
\label{1.8}
0&\leq \alpha < 2, \\
\label{1.9} b_k&=\bar{c_k}=k^{\alpha} .
\end{align}
If $\alpha=0$ in (\ref{1.9}), we have the Mathieu matrices. They
arise if Fourier's method is used to analyze the Hill--Mathieu
operator on $I=[0,\pi]$
$$
Ly =-y^{\prime \prime} + 2a (\cos 2x) y , \quad y(\pi)= y(0), \; \;
y^\prime(\pi)= y^\prime(0).
$$
In this case J. Meixner and F. W. Sch\"afke proved  (\cite{MS54}, Thm
8, Section 1.5; \cite{MSW80}, p. 87) the inequality $ R_n \leq C n^2
$ and conjectured that the asymptotic $R_n \asymp n^2$ holds. This
has been proved 40 years later by H. Volkmer \cite{V96}.

But what can be said if $0 < \alpha <2? $ Proposition 4 in
\cite{DM14}  shows that if (\ref{1.1})--(\ref{1.3}) and (\ref{1.8})
hold, then
\begin{equation}
\label{1.1b} R_n \geq c n^{1-\alpha}.
\end{equation}

 This estimate from below
cannot be improved in the class (\ref{1.1})--(\ref{1.3}),
(\ref{1.8}) as examples in Section 4 show. But in the special case
(\ref{1.8})--(\ref{1.9}) one could expect the asymptotic
\begin{equation}
\label{1r} R_n \asymp n^{2-\alpha}.
\end{equation}
 We show that
$$
 R_n \leq C n^{2-\alpha},
$$
at least for $0 < \alpha < 11/6.$

Notice that in the Hill--Mathieu case we have $\alpha =0, \,b_k =1 \;
\forall k,$ so the operator $B$ is bounded, while it could be
unbounded in the case $\alpha >0.$ We use the approach of Meixner and
Sch\"afke \cite{MS54}, but complement it with an additional argument
to help us deal with the cases where the operator $B$ is unbounded
(but relatively compact with respect to $L$). The main result is the
following.

\begin{thm}
\label{maintheorem} If the conditions (\ref{1.2}) and (\ref{1.9})
hold, then for each $\alpha \in [0,\frac{11}{6})$ there exist
constants $C_\alpha >0  $ and $N_\alpha \in \mathbb{N} $ such that
\begin{align}
\label{1.12} R_n \leq C_\alpha n^{2-\alpha},   \quad   n \geq
N_\alpha.
\end{align}
\end{thm}
Proof is given in Section \ref{ANUPPERBOUND}. It has two parts. In
Section \ref{UPPERBOUNDAKSECTION},  we prove an upper bound for
Taylor coefficients $|a_k(n)|$ in terms of $k, \,n, \, R_n$ and
$\alpha$ (see Theorem \ref{akbound}). In Section \ref{ANUPPERBOUND}
we show how a certain lower bound on $|a_k(n)|$ , in terms of $k, n$,
and $\alpha$, can be used to prove the desired inequality on
particular subsets of $[0,2)$. In the same section we provide such
lower bounds for $|a_2(n)|, |a_4(n)|, \ldots, |a_{12}(n)|$. This
general scheme could be used in an attempt to prove (\ref{1.12}) for
larger subsets of $[0,2)$. One would then need to compute (and
manipulate) $a_k(n)$ for values of $k>12$. See Section
\ref{ANUPPERBOUND} for details.

\section{An upper bound for $|a_k(n)|$}
\label{UPPERBOUNDAKSECTION}

In what follows in this section, suppose that $n$ is a {\em fixed}
positive integer.
\begin{thm}
\label{akbound}  In the above notations, and under the conditions
(\ref{1.2}) and (\ref{1.3}), if

(a)   $ \alpha \in [0,2) $ and (\ref{1.5}) holds,  or \hspace{3mm}
(b) $ \alpha \in [0,1), $ \\ then
\begin{align}
\label{2.1} |a_k(n)| &\leq C \rho^{-(k-1)} \left( n^\alpha +
\rho^{\frac{\alpha}{2-\alpha}} \right), \quad  0 < \rho < R_n,
\end{align}
where $C=C(\alpha,M).$
\end{thm}

\begin{proof}
For $r>0,$ let $$\Delta_r=\{z \in \mathbb{C}: |z| < r \}, \quad
C_r=\{z\in \mathbb{C}:|z|=r\}.$$

Let us choose, for every $z \in \Delta_{R_n},$ an eigenvector $g(z)
= (g_n (z))_{n=1}^\infty $ such that $ \|g(z)\|_{\ell^2} =1 $ (this
is possible by Proposition \ref{prop01}). Then
\begin{equation}
\label{2.2} (L+zB)g(z)=E_n (z) g(z),     \quad  \|g(z)\|_{\ell^2}
=1,
\end{equation}
which implies (after multiplication from the right by $g(z)$)
\begin{equation}
\label{2.3} \ell (z) + z b(z) = E_n (z), \quad   z \in \Delta_{R_n},
\end{equation}
where
\begin{equation}
\label{2.7}  \ell(z) : =\langle Lg(z),g(z)
\rangle=\sum_{k=1}^{\infty} k^2 |g_k(z)|^2,
\end{equation}
and
\begin{equation}
\label{2.8}  b(z) : =\langle Bg(z),g(z) \rangle=\sum_{k=1}^{\infty}
 \left ( c_k g_k (z) \overline{g_{k+1} (z)} + b_k g_{k+1} (z)
\overline{g_k (z)}\right ).
\end{equation}

The functions $\ell (z) $ and $b(z)$  are bounded if $|z|\leq \rho
<R_n.$ Indeed, by (\ref{2.7}) we have $\ell (z)>0.$  By (\ref{2.8})
and (\ref{1.3})
\begin{equation}
\label{2.13} |b(z)| \leq \sum_{k=1}^\infty Mk^\alpha \left(|g_k
(z)|^2 + |g_{k+1}|^2\right) \leq 2M \sum_{k=1}^\infty k^\alpha |g_k
(z)|^2,
\end{equation}
so, estimating the latter sum by H\"older's inequality, we get
\begin{equation}
\label{2.14} |b(z)| \leq 2M (\ell(z))^{\alpha/2}.
\end{equation}
Therefore, in view of (\ref{2.3}).
$$
\ell (z) \leq |E_n (z) | + |z b(z)| \leq |E_n (z) |+ 2M \rho
(\ell(z))^{\alpha/2},   \quad |z| \leq \rho.
$$
Now, Young's inequality implies
$$
\ell (z) \leq |E_n (z) | + (1-
\alpha/2)2^{\frac{\alpha}{2-\alpha}}(2M\rho)^{\frac{2}{2-\alpha}}+
(\alpha/4)\cdot \ell(z),
$$
so, in view of (\ref{1.8}), $\ell (z) $ is bounded by
$$
\ell (z) \leq  2|E_n (z) | + 2(1-
\alpha/2)2^{\frac{\alpha}{2-\alpha}}(2M\rho)^{\frac{2}{2-\alpha}},
\quad |z| \leq \rho.
$$
By (\ref{2.14}), the function $b(z) $ is also bounded if $|z| \leq
\rho.$

Since in (\ref{2.2}) the vectors $g(z), \;z \in \Delta_{R_n}, $  are
chosen in an arbitrary way, we cannot expect the function $z \to
g(z)$ to be continuous, or even measurable. But the functions $\ell
(z)$ and $b(z)$ are measurable. The explanation of this fact is the
only difference in the proof of  (\ref{2.1}) in the cases (a) and
(b).

(a) The functions $\ell (z)$ and $b(z)$  are continuous on
$\Delta_{R_n}\setminus (-R_n, R_n).$

Indeed, in view of (\ref{2.8}) the symmetry assumption (\ref{1.5})
implies that the function $b(z)$ is real--valued.  Therefore, from
(\ref{2.3}) it follows $y b(z) = Im \, E_n (z)$  with $z=x+iy,$ so
$\ell (z)$ and $b(z)$ are continuous on $\Delta_{R_n}\setminus (-R_n,
R_n)$ because
\begin{equation}
\label{2.9}
 b(z) = \frac{1}{y} \, Im ( E_n (z)), \quad \ell (z) = Re (E_n (z)) -
  \frac{x}{y} \, Im ( E_n (z)), \quad   y\neq 0.
\end{equation}

(b) For every $z$ such that $E_n (z) $ is a simple eigenvalue of $L+
wB$ the values  $\ell (z)$ and $b(z)$ are uniquely determined by
(\ref{2.7}) and (\ref{2.8}) and do not depend on the choice of the
vector $g(z) $ in (\ref{2.2}). Therefore, the functions $\ell (z)$
and $b(z)$ are uniquely determined on the set
$$
U= \{z \in \Delta_{R_n}: \;\; E_n (z) \;\; \text{is a simple
eigenvalue of} \;\;L+ zB \}.
$$
On the other hand, the set $ \Delta_{R_n} \setminus U$ is at most
countable and has no finite accumulation points (see Section 5.1 in
\cite{DM14}).

If $w \in U,$ then it is known (\cite{Kato}, Ch.VII, Sect. 1-3, in
particular, Theorem 1.7) that there is a disc $D(w,\tau) $ with
center $w$ and radius $\tau$ such that $E_n (z) $ is a simple
eigenvalue of the operator $L+ zB $ for $z\in D(w,\tau) $ and there
exists an analytic eigenvector function $\psi (z) $ defined in
$D(w,\tau), $  i.e.,
$$
(L+zB)\psi (z) = E_n (z) \psi (z),  \quad \psi (z) \neq 0, \quad z
\in D(w,\tau).
$$
Let $g(z)= \psi (z)/\|\psi (z)\|_{\ell^2}  $  for $ z \in
D(w,\tau).$ Then the coordinate functions $g_k (z) $ are continuous,
and by (\ref{2.7}) the function $\ell (z), \; z \in D(w,\tau),$ is a
sum of a series of positive continuous terms. Therefore, the
function $\ell (z)$ is lower semi--continuous in $D(w,\tau),$ so it
is lower semi--continuous in $U.$ Thus, $\ell (z)$ is measurable on
$\Delta_{R_n}.$ By (\ref{2.3}) we have $b(z) = (E_n (z) - \ell
(z))/z $ for $z\neq 0.$  Thus, $b(z)$ is  measurable in
$\Delta_{R_n}$  as well.

For each $\rho \in (0,R_n),$ consider the space $L^2 (C_\rho)$ with
the norm $\|\cdot\|_\rho$ defined by $ \|f\|_\rho^2 = \frac{1}{2\pi}
\int_{0}^{2 \pi} |f(\rho e^{i \theta})|^2 d\theta. $ The functions
$\ell (z)$ and $b(z)$ are integrable on each circle $C_\rho, \; \rho
<R_n $ because they are bounded and measurable on $C_\rho.$

From (\ref{2.14}) and H\"older's inequality it follows that
\begin{equation}
\label{2.15} \|b(z)\|_\rho \leq 2M \|\ell(z)\|_\rho^{\alpha/2}.
\end{equation}

Since $\ell (z) >0, $ by (\ref{2.3}) and (\ref{2.14}) we have
$$
|Im \, (E_n (z) - n^2)|=|Im \, (z b(z))| \leq
 \rho | b(z)|.
$$
Therefore,
\begin{equation}
\label{2.16} \|Im \,(E_n (z) - n^2)\|_\rho \leq \rho \cdot \|
b(z)\|_\rho.
\end{equation}
If $f$ is an analytic function defined on $\Delta_{R_n}$ with
$f(0)=0$, then $\|Re(f)\|_\rho=\|Im(f)\|_\rho$. In particular, we
have
$$
\|Re \,(E_n (z) - n^2)\|_\rho = \|Im \,(E_n (z) - n^2)\|_\rho,
$$
which implies, by (\ref{2.16}),
\begin{equation}
\label{2.18}
 \|E_n (z) -n^2\|_\rho \leq \sqrt{2} \rho \cdot  \| b(z)\|_\rho.
\end{equation}
In view of (\ref{2.3}) and (\ref{2.18}), the triangle inequality
implies
$$
\|\ell\|_\rho \leq n^2 + \|E_n (z) -n^2 \|_\rho + \|b(z)\|_\rho \leq
n^2 + (1+\sqrt{2})\rho \cdot  \| b(z)\|_\rho.
$$
Therefore, from (\ref{2.15}) it follows that
\begin{equation}
\label{2.21} \|\ell\|_\rho \leq  n^2 +5 M \rho \| \ell
\|_\rho^{\alpha/2}.
\end{equation}
Now, Young's inequality yields
$$
5 M \rho \| \ell \|_\rho^{\alpha/2} \leq \left
(1-\alpha/2)(5M2^{\alpha/2} \rho \right)^{\frac{2}{2-\alpha}} +
\frac{\alpha}{4} \|\ell\|_\rho \leq C_1 \rho^{\frac{2}{2-\alpha}} +
\frac{1}{2}\|\ell\|_\rho,
$$
with $C_1 = \left (1-\alpha/2)(5M2^{\alpha/2}
\right)^{\frac{2}{2-\alpha}}.$  Thus, by (\ref{2.21}), we have
$$
 \|\ell\| \leq 2 n^2 + 2C_1 \rho^{\frac{2}{2-\alpha}}.
$$
In view of (\ref{2.18}) and (\ref{2.15}), this implies
\begin{equation}
\label{2.22} \|E_n (z) -n^2\|_\rho \leq 3 M \rho \left (2^{\alpha/2}
n^{\alpha} + (2C_1)^{\alpha/2} \rho^{\frac{\alpha}{2-\alpha}}
\right).
\end{equation}

By Cauchy's formula, we have
$$ a_k(n)=\frac{1}{2 \pi i} \int_{\partial \Delta_\rho}
\frac{E_n(\zeta)-n^2}{\zeta^{k+1}} d\zeta.
$$
From  (\ref{2.22}) it follows that
$$
|a_k(n)| \leq \rho^{-k}\|E_n (z) -n^2\|_\rho
         \leq 3 M  \rho^{-k+1}  \left (2^{\alpha/2}
n^{\alpha} + (2C_1)^{\alpha/2} \rho^{\frac{\alpha}{2-\alpha}}
\right),
$$
which implies (\ref{2.1}) with $C= 3M (2+ 2C_1)^{\alpha/2}.$ This
completes the proof of Theorem~\ref{akbound}.
\end{proof}

{\em Remark.}  In fact, to carry out the proof of
Theorem~\ref{akbound} we need only to know that there exists a pair
of functions $\ell (z)$ and $b(z)$ which satisfy (\ref{2.3}) and
(\ref{2.14}), and are integrable on each circle $C_\rho, \; \rho
<R_n. $ We explained that the pair defined by (\ref{2.2}),
(\ref{2.7}) and (\ref{2.8}) has these properties. In the case (a) of
Theorem~\ref{akbound} the same argument could be used to define a
pair of real analytic functions functions $\ell (z)$ and $b(z)$ which
satisfy (\ref{2.3}) and (\ref{2.14}).

Indeed, by (\ref{1.5}) the operator $B$ is a self--adjoint, so $L+ x
B, \; x \in \mathbb{R},$ is self--adjoint as well. Thus, the function
$E_n (z)$ takes real values on the real line and its Taylor's
coefficients are real. Since the quotients $ \frac{1}{y} \,Im
(x+iy)^k, \; k \in \mathbb{N}, $ are polynomials of $y,$ it is easy
to see by the Taylor series of $E_n (z)$ that $\frac{1}{y} Im ( E_n
(z))$ (defined properly for $y=0$) is a real analytic function in
$\Delta_{R_n}.$  Therefore, if one defines a pair of functions
$\tilde{\ell} (z)$ and $\tilde{b}(z)$ by (\ref{2.9}), then
(\ref{2.3}) holds immediately, and (\ref{2.14}) follows because on
$\Delta_{R_n}\setminus (-R_n, R_n)$ these functions coincide with
$\ell (z)$ and $b(z).$

\section{An upper bound for $R_n$} \label{ANUPPERBOUND}
In this section we use (\ref{2.1}) in the case of (\ref{1.9}) to
prove Theorem \ref{maintheorem}. Roughly speaking, the bound
(\ref{1.12}) will be achieved for $\alpha \in [0,\frac{11}{6})$ by
inserting the known (from \cite{DM14}) formulas for $a_2(\alpha,n),
\ldots, a_{12}(\alpha,n)$ into inequality (\ref{2.1}). With our
approach, using only $a_{2k}, \; k \leq 6,$ it is possible to get
good lower bounds only if $0 \leq \alpha < 11/6.$

We begin with the following observation.
\begin{lem}
\label{lem33} Suppose the conditions (\ref{1.2}),(\ref{1.3}) and
(\ref{1.8}) hold.

 (a)    If for some fixed $k,n \in \mathbb{N}$  and
$\alpha \in [0,2-\frac{2}{k})$ we have
 $a_k (n) \neq 0, $ then $R_n < \infty.$

(b)  If $R_n =\infty,$ then $E_n (z) $ is a polynomial such that
$\deg E_n (z) \leq \frac{\alpha}{2-\alpha}.$
\end{lem}

\begin{proof}
Let $a=|a_k (n) |>0.$  Then, by Theorem~\ref{akbound},
\begin{equation}
\label{3.01} a \rho^{k-1} \leq C \left (   n^\alpha +
\rho^{\frac{\alpha}{2-\alpha}}  \right ), \quad   \forall \, \rho<
R_n.
\end{equation}
The condition $\alpha \in [0,2-\frac{2}{k})$ implies  $k-1 >
\frac{\alpha}{2-\alpha};$ therefore, (\ref{3.01}) fails for
sufficiently large $\rho.$ Thus, $R_n \leq \sup\{\rho:  \; \rho \in
(\ref{3.01}) \} < \infty,$ which proves (a).

If $R_n = \infty,$  then (a) shows that $a_k (n) = 0 $ for all $k$
such that $k > \frac{\alpha}{2-\alpha}.$  This proves (b).

\end{proof}

\begin{lem}
\label{lem3} Suppose that conditions (\ref{1.2}) and (\ref{1.3})
hold. If for some fixed $k,n \in \mathbb{N},$ $A>0$ and $\alpha \in
[0,2-\frac{2}{k})$ we have
\begin{align}
\label{3.1a}
 A n^{k \alpha - 2(k-1)} \leq |a_{k}(n)|,
\end{align}
then
\begin{align}
\label{3.1} R_n  \leq  \tilde{C} n^{2-\alpha},
\end{align}
where $\tilde{C}=\tilde{C} (\alpha, M, A,k).$

\begin{proof}
It is enough to prove that
\begin{equation}
\label{3.2} \rho  \leq  \tilde{C} n^{2-\alpha}, \quad \forall \,
\rho \in (0, R_n).
\end{equation}
Then (\ref{3.1}) follows if we let $\rho \to R_n.$

By (\ref{2.1}) we have
$$
A n^{k \alpha - 2(k-1)} \leq | a_{k}(n) |
                     \leq 2 C(\alpha,M) \rho^{-(k-1)}
\max(n^{\alpha},\rho^{\frac{\alpha}{2-\alpha}}).
$$
If $n^{\alpha} \geq  \rho^{\frac{\alpha}{2-\alpha}}$, then we get
(\ref{3.2}) with $\tilde{C}=1.$

Suppose that $n^{\alpha} <  \rho^{\frac{\alpha}{2-\alpha}}.$ Then
$\max(n^{\alpha},\rho^{\frac{\alpha}{2-\alpha}})=
\rho^{\frac{\alpha}{2-\alpha}},$
 so
$$
A \rho^{k-\frac{2}{2-\alpha}} \leq 2C(\alpha,M)
(n^{2-\alpha})^{k-\frac{2}{2-\alpha}}.
$$
Thus, whenever $\alpha < 2 - 2/k,$ this inequality implies
(\ref{3.1}) with $\tilde{C} = (2C/A)^\gamma,$ where $ \gamma =
(2-\alpha)/(k(2-\alpha)-2).$
\end{proof}
\end{lem}

According to the preceding lemma, all one needs in order to get an
upper bound on $R_n$ of the form (\ref{3.1}) (or even to explain
that $R_n$ is finite) is to find a lower bound on $|a_k (n)|$ of the
form (\ref{3.1a}) (or at least to explain that $a_k (n) \neq 0). $
We now describe a technique to provide such lower bounds. Theorem
\ref{maintheorem} will follow when we get such lower bounds for
$|a_2 (n)|, \ldots, |a_{12} (n)|$.

\begin{lem}
Under conditions (\ref{1.4}) and (\ref{1.9}), for each fixed
$\alpha<2,$ the coefficient $a_k(n,\alpha)$ can be written in the
form
\begin{align}
\label{3.8} a_k(n,\alpha)=n^{k\alpha-(k-1)}f_{\alpha}(1/n)
\end{align}
where
$$
f_{\alpha}(w)=\sum_{j=0}^{\infty}P_k(j,\alpha) w^{j}
$$
is analytic on the disk $ |w| < 1/k,$  and $P_k(j,\alpha)$ are
polynomials of $\alpha. $

\begin{proof}
We begin this proof by stating the equation (3.7) from \cite{DM14}
\begin{align}
\label{3.10} a_k(n) =  \frac{1}{2\pi i}
          \int_{\partial \Pi} \left( \sum_{|j-n| \leq k} (\lambda - n^2)
           \langle R_{\lambda}^0 (B R_{\lambda}^0)^ke_j,e_j \rangle
           \right) d\lambda,
\end{align}
where $R_{\lambda}^0= (\lambda - L)^{-1}$, $e_j$ is the $j^{th}$
unit vector, and $\Pi$ is the square centered at $n^2$ of width
$2n$. This formula appears in \cite{DM14} only in the case of
$\alpha \in [0,1)$, but its proof therein holds for $\alpha < 2$ as
well. It follows from (\ref{1.1}) that for each $j \in N$,
$$
 BR_{\lambda}^0 e_j = \begin{cases}\frac{(j-1)^\alpha
}{\lambda - j^2} e_{j-1} +
\frac{j^\alpha }{\lambda - j^2} e_{j+1} &\text{ if j $>$ 1} \\ \\
\frac{ 1}{\lambda - 1} e_2  &\text{ if j $=$ 1}.
\end{cases}
$$
So, $(\lambda-n^2) \langle R^0_\lambda(BR^0_\lambda)^k e_j,e_j
\rangle$ can be written as a finite sum each of whose terms is of the
form
$$\frac{\lambda-n^2}{\lambda - (n-j^\prime_0)^2}
\prod_{i=1}^{k}\frac{(n-d^\prime_i)^{\alpha}}{\lambda-(n-j^\prime_i)^2}
$$
with $j^\prime_i$ and $d^\prime_i$ integers satisfying $|j^\prime_i|,
|d^\prime_i| < k$ for each $i$. So, from a residue calculation on
(\ref{3.10}), $a_k (n)$ can be written as a linear combination of
terms of the form
\begin{equation}
\label{3.12} (n-d_k) ^{\alpha}\prod_{i=1}^{k-1}
\frac{(n-d_i)^\alpha}{n^2 - (n-j_i)^2}
\end{equation}
$$
=Cn^{k\alpha - (k-1)}\left(1-\frac{d_k}{n} \right)^\alpha
\prod_{i=1}^{k-1} \left[ \left(1-\frac{d_i}{n}\right)^\alpha
\left(1-\frac{j_i}{2n} \right)^{-1} \right]
$$
with $C= \prod_{i=1}^{k-1}  (2j_i)^{-1} $ and $|j_i|, |d_i| < k$ for
each $i$.

For $n > k$, we have $|d_i/ n| <1$ and $|j_i/(2 n)| < 1.$  Thus,
\begin{align}
\label{3.14} \left(1-\frac{d_i}{n}\right)^\alpha &= 1-\alpha
\left(\frac{d_i}{n}\right) +
 \frac{\alpha(\alpha-1)}{2} \left( \frac{d_i}{n}\right)^2 + \ldots
\\
\label{3.15} \left(1-\frac{j_i}{2n}\right)^{-1}   &= 1 + \left(
\frac{j_i}{2n}\right) +
 \left( \frac{j_i}{2n}\right)^2 + \ldots
\end{align}
 are analytic functions of $z=1/n$ whenever $n > k$.
 Combining (\ref{3.12}) with (\ref{3.14})--(\ref{3.15}),
 we deduce that $a_k (n)$ can
 be written as in (\ref{3.8})
 with $f_\alpha (z) $ analytic for $|z|< 1/k$.
\end{proof}
\end{lem}

The preceding lemma guarantees that whenever $\alpha <2$,
$$
a_k(n,\alpha) = P_k(0,\alpha) n^{k\alpha-(k-1)} + O(n^{k\alpha-k})
\quad \text{as} \,\,\, n \rightarrow \infty.
$$
When $a_2 (n),  \ldots, a_{12}(n)$ were computed (following the
approach of \cite[p.305--306]{DM14}), an interesting phenomenon was
observed. If $2 \leq k \leq 12$, then
\begin{align}
P_k(j,\alpha)=0 \,\,\,  \text{for each} \,\,\,0 \leq j \leq k-2.
\end{align}
In particular, if (\ref{1.8}) and (\ref{1.9}) hold, then
\begin{align}
\label{3.18} a_{k}(n) &= P_{k}(k-1,\alpha)n^{k\alpha-2(k-1)} +
O(n^{k\alpha-2k+1}), \,\,\,\, n\rightarrow \infty;
\end{align}
the polynomials $P_{k}(k-1,\alpha), \; k=2, 4,\ldots, 12, $ are given
in the following table. \vspace{15mm}

\begin{tabular}{r|l}
k       &$P_{k}(k-1,\alpha)$\\
\hline \\
$2$& $-\alpha + \frac{1}{2}$\\ \\
\hline \\
$4$& $-\alpha^3+\frac{9}{4}\alpha^2-\frac{11}{8}\alpha+\frac{5}{32}$\\ \\
\hline \\
$6$& $-\frac{9}{4}\alpha^5+\frac{73}{8}\alpha^4-\frac{27}{2}\alpha^3+\frac{281}{32}\alpha^2-\frac{147}{64}\alpha+\frac{9}{64}$\\ \\
\hline \\
$8$& $-\frac{61}{9} \alpha^7 + \frac{2881}{72}\alpha^6
           -\frac{6875}{72}\alpha^5 + \frac{33937}{288}\alpha^4 - \frac{11437}{144} \alpha^3
        + \frac{64649}{2304}\alpha^2 -\frac{4507}{1024} \alpha +\frac{1469}{8192} $ \\ \\
\hline \\
$10$&$-\frac{1525}{64} \alpha^9 + \frac{23705}{128} \alpha^8-
\frac{353023}{576} \alpha^7
         + \frac{648539}{576} \alpha^6 - \frac{5774039}{4608} \alpha^5 + \frac{7955297}{9216} \alpha^4$ \\ \\&

         $- \frac{6626165}{18432} \alpha^3 + \frac{6173425}{73728} \alpha^2- \frac{148881}{16384} \alpha
         + \frac{4471}{16384}$ \\ \\
\hline \\
$12$&$ - \frac{221321}{2400}\alpha^{11} +
\frac{8544347}{9600}\alpha^{10}
          -\frac{1207947}{320}\alpha^9 +\frac{71029219}{7680}\alpha^8
          - \frac{92577243}{6400} \alpha^7 + \frac{385333821}{25600}\alpha^6$ \\ \\&
         $- \frac{16162765}{1536} \alpha^5 + \frac{9344339}{1920} \alpha^4
          -\frac{583689039}{409600} \alpha^3 + \frac{296768801}{1228800}\alpha^2
          - \frac{12877899}{655360} \alpha  + \frac{121191}{262144}$
\end{tabular}
\\ \\ \\ \\

Numerical computations tell us that in the following table, each
inequality in the second column holds on the union of intervals shown
in the first column.
\\ \\
\begin{tabular}{l|l}
\\
Set & Inequality
\\
\hline
\\
$\alpha \in S_2 = \left[ 0,\frac{1}{4} \right] \cup \left[
\frac{3}{4},1 \right)$
& $|P_2(1,\alpha)|> \frac{1}{8}$\\ \\
\hline \\
$\alpha \in S_4 = \left[ \frac{1}{4},\frac{3}{4} \right] \cup \left[
1,\frac{9}{8} \right] \cup \left[ \frac{11}{8},\frac{3}{2}
\right)$ & $|P_4(3,\alpha)|>\frac{1}{32}$\\  \\
\hline \\
$\alpha \in S_6 = \left[ \frac{9}{8},\frac{11}{8} \right] \cup \left[
\frac{25}{16},\frac{5}{3} \right)$ &
$|P_6(5,\alpha)| > \frac{1}{200}$\\  \\
\hline \\
$\alpha \in S_8 = \left[ \frac{3}{2},\frac{25}{16} \right] \cup
\left[ \frac{5}{3},\frac{7}{4} \right)$ &
$|P_8(7,\alpha)| > \frac{1}{10}$\\  \\
\hline \\
$\alpha \in S_{10} = \left[ \frac{7}{4},\frac{9}{5}
\right)$ & $|P_{10}(9,\alpha)|>\frac{1}{2}$\\  \\
\hline  \\
$\alpha \in S_{12} = \left[ \frac{9}{5},\frac{11}{6} \right) $ & $|P_{12}(11,\alpha)| > 1$\\  \\
\hline \\
\end{tabular}
\\ \\ \\ \\

{\em Proof of Theorem \ref{maintheorem}.} In view of (\ref{3.18}) and
the above table, there is a constant $A>0$  such that, for each
$\alpha \in [0,2-\frac{1}{6}),$ we have
\begin{equation}
\label{3.30} |a_k(n,\alpha)|>A n^{k\alpha-2(k-1)}, \quad n \geq
N_\alpha.
\end{equation}

Therefore, Lemma \ref{lem3} implies that there exists a constant
$C_{\alpha}$ such that $$R_n \leq C_\alpha n^{2-\alpha} \quad
\text{for}  \; n \geq N_\alpha.$$ Thus,
 (\ref{1.12}) holds for  $n\in \mathbb{N},$
which completes the proof of Theorem~\ref{maintheorem}.

\section{General discussion}
In this section we give a few examples to show that the order $1-
\alpha $ of lower bound (\ref{1.1b}) for $R_n$ is sharp in the class
of matrices $B$ with (\ref{1.2})--(\ref{1.4}).

1. {\it A case in which $R_n \sim n^{1-\alpha}.$} Let  $\alpha \in
[0,1).$ Suppose now that in (\ref{1.1}) we set
\begin{align}
\label{4.1}
b_k&=c_k=(2+(-1)^k)k^{\alpha} \\
\label{4.2} q_k&=k^2
\end{align}
Then by \cite{DM14}, Section 7.5, p.35,
$$
|a_2(n)| = \left| \frac{b_{n-1} c_{n-1}}{2n-1} - \frac{b_n c_n}{2n+1}
\right|
       = \begin{cases} \left| \frac{9(n-1)^{2\alpha}}{2n-1}-
       \frac{n^{2\alpha}}{2n+1} \right| & \text{if $n$ is odd,} \vspace{1mm}\\
                         \left| \frac{(n-1)^{2\alpha}}{2n-1}-
                         \frac{9 n^{2\alpha}}{2n+1} \right| & \text{if $n$ is even}
                         \end{cases}
$$
so
$$
|a_2(n)|\geq c \,n^{2\alpha - 1}, \quad c>0.  $$ In view of  Lemma
\ref{lem33}, this implies that $R_n < \infty $ for $\alpha \in
[0,1).$

 Therefore, by (\ref{2.1}) in Theorem \ref{akbound}, for
each $\alpha \in [0,1),$ we have
\begin{equation}
\label{5.1} n^{2\alpha - 1} \leq |a_2(n)| \leq  2C(\alpha) R_n^{-1}
\max(n^{\alpha},R_n^{\frac{\alpha}{2-\alpha}}), \quad n \geq n_0.
\end{equation}

If $n^\alpha \leq R_n^{\frac{\alpha}{2-\alpha}},$ then $R_n \geq
n^{2-\alpha} $ and (\ref{5.1}) gives $ n^{2\alpha-1} \leq 2C(\alpha)
R_n^{\frac{2\alpha-2}{2-\alpha}}, $ which implies
$$
2C(\alpha) \geq n^{2\alpha-1} R_n^{\frac{2-2\alpha}{2-\alpha}} \geq
n^{2\alpha-1}n^{2-2\alpha}=n.
$$
Therefore, we have $ \max(n^{\alpha},R_n^{\frac{\alpha}{2-\alpha}})=
n^{\alpha} \quad \text{for} \quad n>2C(\alpha). $ So, (\ref{5.1})
implies
\begin{equation}
\nonumber \label{5.2} R_n \leq 2C(\alpha) n^{1-\alpha} \quad
\text{for} \quad n>2C(\alpha).
\end{equation}

 On the other hand, by Proposition 4 of \cite[p.296]{DM14},
we have $ R_n \geq \frac{1}{8} n^{1-\alpha} $ for large enough $n.$
Hence, we have shown that in the special case of
(\ref{4.1})--(\ref{4.2}),
\begin{align}
\label{4.11} R_n \asymp n ^{1-\alpha}.
\end{align}

2. Of course we can simplify the example (\ref{4.1}) by choosing
\begin{equation}
\label{4.12} b_k = c_k = \left [ 1+ (-1)^{k-1}   \right ] k^\alpha
\end{equation}
This ensures that $L+zB - E(z) I$ has the structure of a
tri--diagonal matrix with $2\times2$ blocks along the diagonal. The
$m^{\text{th}}$ block will have the form
\begin{align}
\label{4.13}
\begin{bmatrix}
  T-E     &   z b    \\
  z b   &   V-E      \\
\end{bmatrix},
\end{align}
where
$$
T =(2m-1)^2, \quad V =(2m)^2,  \quad  b = (2m-1)^{\alpha}, \quad
m=1,2, \ldots.
$$
It follows that the two eigenvalues corresponding to this block are
$$
E(z)=\frac{1}{2} \left( T+V \pm \sqrt{(T-V)^2+4 z^2 b^2}\right).
$$
So, the branching points of these branches of $E(z)$ occur at
\begin{equation}
\label{4.14} z_{1,2}=\pm i \left( \frac{V-T}{2b} \right).
\end{equation}
Hence, we have
\begin{equation}
\label{4.15}
z^m_{1,2} = \pm   \frac{i(4m-1)}{2(2m-1)^\alpha} \\
          =\pm i  (2m)^{1-\alpha}
          \left( 1+ \frac{2\alpha-1}{4m}
          +O(m^{-2}) \right)
\end{equation}
Therefore,
$$
R_{2m-1} =R_{2m} \sim (2m)^{1-\alpha},
$$
i.e., we have the same sharp order  $1-\alpha $ as in (\ref{4.11}).
\vspace{3mm}

3. This simplified example (\ref{4.12}) is extreme in the sense that
the spectral Riemann surface (SRS)
$$
G(B) = \{(z,E) \in \mathbb{C}^2 : \quad (L+zB)f = Ef, \quad f \in
\ell^2, \; f\neq 0 \}
$$
splits: it is a union of Riemann surfaces defined by determinants of
the blocks (\ref{4.13}), i.e.,
$$
E^2 - E [(2m-1)^2 + (2m)^2] + (2m-1)^2  (2m)^2 -z^2 (2m-1)^2=0, \quad
m\in \mathbb{N}.
$$
In the case (\ref{4.1})  we have no elementary reason to say
anything about (ir)reducibility of the spectral Riemann surface
$G(B)\;$ (see more about irreducibility of SRS in \cite{DM14,V04}).

Nevertheless, we would conjecture that this surface $G(B)$ is {\em
irreducible} if $B\in (\ref{4.1}),$ or more generally, if
\begin{equation}
b_k = c_k \left ( 1+\gamma (-1)^{k-1} \right ) k^\alpha, \quad 0\leq
\gamma <1.
\end{equation}
If $\gamma=0 $ we proved  in \cite{DM14}, Theorem 3,  such
irreducibility for $\alpha=1/2$ and many but not all $\alpha's $ in
$[0;1/2].$ \vspace{3mm}

If $1 \leq \alpha < 2$ let us choose in (\ref{4.13})
\begin{equation}
b = b_m = \frac{1}{B_m}  (2m-1)^{\alpha}, \quad |B_m|\geq 1.
\end{equation}
Then (\ref{4.14}) holds, so by (\ref{4.15})
$$
z_{1,2} = \pm i B_m (2m)^{1-\alpha}\left( 1 + O(1/m) \right).
$$
The sequence $\{B_m\}$ could be chosen in such a way that the set $A$
of accumulation points for $\{z_{1,2}^m\}$ is the entire complex
plane $\mathbb{C}$, or for any closed $K \subset \mathbb{C}$ with
$K=-K$ we can make $A=K.$

4. Our argument in Section 2, uses Young's and H\"older's
inequalities, i.e., the concavity of the function $x^{\alpha/2}, \; 1
\leq x < \infty, \; 0 \leq \alpha  < 2.$ It cannot be applied if
$\alpha < 0$ although in this case the operator $B\in (\ref{1.3})$ is
even compact. {\em Yet, we conjecture that $R_n \leq K(\alpha)
n^{2-\alpha}$ holds both for $\alpha \in [\frac{11}{6},\, 2)$ and
$\alpha <0.$ Moreover, we expect that our conjecture (\ref{1r}) holds
for $\alpha <0$  as well.}

\end{document}